\documentclass[10pt,draft]{amsart}
\usepackage{amsfonts}
\usepackage{ifthen}
\usepackage{amsthm}
\usepackage{amsmath}
\usepackage{graphicx}
\usepackage{amscd,amssymb,amsthm}
\usepackage{enumerate}
\usepackage{color}

\newcounter{minutes}
\divide\time by 60
\newcounter{hours}
\multiply\time by 60 \addtocounter{minutes}{-\time}

\title[Tur\'an type inequalities for confluent hypergeometric functions]{Tur\'an type inequalities for confluent hypergeometric functions of the second kind}

\author[\'A. Baricz]{\'Arp\'ad Baricz}
\address{Department of Economics, Babe\c{s}-Bolyai University, Cluj-Napoca 400591, Romania}
\address{Institute of Applied Mathematics, John von Neumann Faculty of Informatics, \'Obuda University, 1034 Budapest, Hungary}
\email{bariczocsi@yahoo.com}

\author[S. Ponnusamy]{Saminathan Ponnusamy}
\address{Indian Statistical Institute, Chennai Centre, Society for Electronic Transactions and Security,
MGR Knowledge City, CIT Campus, Taramani, Chennai 600113, India}
\email{samy@iitm.ac.in}

\author[S. Singh]{Sanjeev Singh}
\address{Department of Mathematics, Indian Institute of Technology Madras, Chennai 600036, India}
\email{sanjeevsinghiitm@gmail.com}

\thanks{$^{\bigstar}$The work of \'A. Baricz was supported by the J\'anos Bolyai Research Scholarship of
the Hungarian Academy of Sciences. The second author is on leave from the Department of Mathematics,
Indian Institute of Technology Madras, Chennai-600 036, India. The research of S. Singh was supported by the
fellowship of the University Grants Commission, India.}

\newtheorem{theorem}{Theorem}

\begin{document}

\def\thefootnote{}
\footnotetext{ \texttt{File:~\jobname .tex,
          printed: \number\year-0\number\month-\number\day,
          \thehours.\ifnum\theminutes<10{0}\fi\theminutes}
} \makeatletter\def\thefootnote{\@arabic\c@footnote}\makeatother


\maketitle


\begin{abstract}
In this paper we deduce some tight Tur\'an type inequalities for Tricomi confluent hypergeometric functions of the second kind, which in some cases improve the existing results in the literature. We also give alternative proofs for some already established Tur\'an type inequalities. Moreover, by using these Tur\'an type inequalities, we deduce some new inequalities for Tricomi confluent hypergeometric functions of the second kind. The key tool in the proof of the Tur\'an type inequalities is an integral representation for a quotient of Tricomi confluent hypergeometric functions, which arises in the study of the infinite divisibility of the Fisher-Snedecor $F$ distribution.
\end{abstract}

\section{\bf Introduction}
\setcounter{equation}{0}
Let us start with the following inequality for the Legendre polynomials
$$P^2_n(x)-P_{n-1}(x) P_{n+1}(x)\geq0,$$
where $x\in[-1,1]$ and $n\in \{1,2\ldots\}$, which was first proved by the Hungarian mathematician Paul Tur\'an \cite{PT}  while studying the zeros of Legendre polynomials. Thereafter an inequality of this kind is known as a Tur\'an type inequality. G. Szeg\H {o} \cite{szego} has given four elegant proofs of the above inequality. There is an immense literature on the Tur\'an type inequalities for  classical (orthogonal) polynomials and special functions. For recent development on the topic we suggest the reader to refer \cite{baricz}-\cite{bps2} and references therein.

In this paper we study the Tur\'an type inequalities for the confluent hypergeometric functions of the second kind denoted by $\psi(a,c,\cdot)$ also known as Tricomi confluent hypergeometric functions. This function is a particular solution of Kummer's differential equation \cite[p. 504]{abramowitz}
$$xy''(x)+(c-x)y'(x)-ay(x)=0,$$
and for $a>0$, $c\in \mathbb{R}$, and $x>0$ has the following integral representation \cite[p. 505]{abramowitz}
\begin{equation}\label{integr}\psi(a,c,x)=\frac{1}{\Gamma(a)}\int_{0}^{\infty}e^{-xt}t^{a-1}(1+t)^{c-a-1}dt.\end{equation}

Recently Baricz and Ismail \cite{bais} have deduced some Tur\'an type inequalities for confluent hypergeometric functions of the second kind. In \cite{baricz4}, the author found some tight bounds for Tur\'anian of modified Bessel functions of first and second kind. Motivated by the results from \cite{baricz4}, in this paper we find tight bounds for the Tur\'anians of confluent hypergeometric functions of the second kind and we offer some alternative proofs of the Tur\'an type inequalities given in \cite{bais}. Moreover, by using a technique similar to \cite{baricz1}, we derive some new inequalities for confluent hypergeometric functions of the second kind. The organization of the paper is as follows: in Sect. \ref{MR} we state our main results, while in Sect. \ref{proof} we prove the main results and we give some alternative proofs of the Tur\'an type inequalities derived in \cite{bais}.

\section{\bf Main Results}\label{MR}
\setcounter{equation}{0}

Our first set of results deals with some tight Tur\'an type inequalities for confluent hypergeometric functions of the second kind, which improve the existing Tur\'an type inequalities proved in \cite{bais}. The first main result is related to the bounds of the Tur\'anian of confluent hypergeometric functions of the second kind when we shift both parameters.

\begin{theorem}\label{theorem2*}
If $a>0$, $c<1$ and $x>0$, then the following Tur\'an type inequality is valid
\begin{equation}\label{te1}
\frac{c-a-1}{x^2}\psi^2(a,c,x)<\psi^2(a,c,x)-\psi(a-1,c-1,x)\psi(a+1,c+1,x).
\end{equation}
Similarly, if $a>1$, $c<-1$ and $x>0$, then we have the next inequality
\begin{equation}\label{te1r}
\psi^2(a,c,x)-\psi(a-1,c-1,x)\psi(a+1,c+1,x)<\left(\frac{1}{c}+2x\left(\frac{c-a}{c^2(c+1)}\right)\right)\psi^2(a,c,x).
\end{equation}
In addition, the following Tur\'an type inequality
\begin{equation}\label{te2}
-\frac{1}{2x}\psi^2(a,c,x)<\psi^2(a,c,x)-\psi(a-1,c-1,x)\psi(a+1,c+1,x)
\end{equation}
holds for all $a>0$, $c<1$ and $x>0$. The inequalities \eqref{te1} and \eqref{te2} are sharp as $x\rightarrow \infty$, while \eqref{te1r} is sharp as $x\rightarrow 0$.
\end{theorem}

We would like to mention here that for $x^2>c(c-a-1)$, the inequality \eqref{te1} is better than the Tur\'an type inequality \cite[Theorem 2]{bais}
\begin{equation}\label{te7}
\frac{1}{c}\psi^2(a,c,x)<\psi^2(a,c,x)-\psi(a-1,c-1,x)\psi(a+1,c+1,x),
\end{equation}
which is valid for $a>0>c$, and $x>0$, while for $x<\frac{c(c+1)}{2(a-c)}$, the inequality \eqref{te1r} is better than the Tur\'an type inequality \cite[Theorem 2]{bais}
\begin{equation}\label{te8}
\psi^2(a,c,x)-\psi(a-1,c-1,x)\psi(a+1,c+1,x)<0
\end{equation}
which is valid for $a>0$, $c<1$ and $x>0$. Moreover, observe that for $x>-\frac{c}{2}>0$, \eqref{te2} is better than \eqref{te7}.

The next result is about the Tur\'anian of confluent hypergeometric functions of the second kind when we shift only the first parameter.

\begin{theorem}\label{theorem3*}
If  $a>0>c$ and $x>0$ then the following Tur\'an type inequalities are valid
\begin{equation}\label{te3}
\frac{1}{1+a-c}\left(1+\frac{1}{2c}x\right)\psi^2(a,c,x)<\psi^2(a,c,x)-\psi(a-1,c,x)\psi(a+1,c,x)
<\frac{2}{x}\psi^2(a,c,x).
\end{equation}
Moreover the right-hand side of \eqref{te3} holds true for all $a>0$, $c<1$ and $x>0$. In addition, the following Tur\'an type inequality
\begin{equation}\label{te5}
\frac{1}{1+a-c}\left(1-\left(\frac{c-a}{c^2(c+1)}\right)x^2\right)\psi^2(a,c,x)
<\psi^2(a,c,x)-\psi(a-1,c,x)\psi(a+1,c,x)
\end{equation}
holds for $a>1$, $c<-1$ and $x>0$. All of the above inequalities are sharp as $x\rightarrow 0$ except right-hand side of \eqref{te3} which is sharp as $x\rightarrow \infty$.
\end{theorem}

We also note that when $x<-\frac{3c}{2}$, the left-hand side of \eqref{te3} improves the Tur\'an type inequality \cite[Theorem 3]{bais}, namely
\begin{equation}\label{te9}
\psi^2(a,c,x)-\psi(a-1,c,x)\psi(a+1,c,x)>0,
\end{equation}
which is valid for $a>0$, $c<1$ and $x>0$, while  for $x>2(1+a-c)$, the right-hand side of \eqref{te3} improves the Tur\'an type inequality \cite[Theorem 3]{bais}
\begin{equation}\label{te10}
\psi^2(a,c,x)-\psi(a-1,c,x)\psi(a+1,c,x)<\frac{1}{1+a-c}\psi^2(a,c,x),
\end{equation}
which is valid $a>1>c$ and $x>0$. Note also that \eqref{te5} is better than \eqref{te9} for $a>1$, $c<-1$ and $x>0$ such that $x^2<\frac{c^2(c+1)}{c-a}$.

Now, we focus on the Tur\'anian of confluent hypergeometric functions of the second kind when we shift only the second parameter.

\begin{theorem}\label{theorem4*}
If $a>0$, $c<1$ and $x>0$ then the following Tur\'an type inequality is valid
\begin{equation}\label{te6}
-\frac{a}{x^2}\psi^2(a,c,x)<\psi^2(a,c,x)-\psi(a,c-1,x)\psi(a,c+1,x).
\end{equation}
In addition, if $a>1$, $c<-1$ and $x>0$ then the following Tur\'an type inequality is valid
\begin{equation}\label{te6r}
\psi^2(a,c,x)-\psi(a,c-1,x)\psi(a,c+1,x)
<\frac{a}{c(1+a-c)}\left(1+2x\left(\frac{c-a}{c(c+1)}\right)\right)\psi^2(a,c,x).
\end{equation}
The inequality \eqref{te6} is sharp as $x\rightarrow \infty,$ while \eqref{te6r} is sharp as $x\rightarrow 0$.
\end{theorem}

For $x^2>c(c-a-1)$, the inequality \eqref{te6} is better than the Tur\'an type inequality \cite[Theorem 4]{bais}
\begin{equation}\label{te11}
\frac{a}{c(1+a-c)}\psi^2(a,c,x)<\psi^2(a,c,x)-\psi(a,c-1,x)\psi(a,c+1,x),
\end{equation}
which is valid for all $a>0>c$ and $x>0$, while for $x<\frac{c(c+1)}{2(a-c)}$, the inequality \eqref{te6r} is better than the Tur\'an type inequality \cite[Theorem 4]{bais}
$$\psi^2(a,c,x)-\psi(a,c-1,x)\psi(a,c+1,x)<0,$$
which is valid for $a,x>0$ and $c\in \mathbb{R}$.

It is worth also to mention that very recently Simon \cite{simon} obtained also some interesting Tur\'an type inequalities for Tricomi confluent hypergeometric functions, however these inequalities cannot be directly comparable with our present results. In \cite{simon} the following Tur\'an type inequalities appear
$$\psi^2(a,c,x)-\psi(a,c-1,x)\psi(a,c+1,x)\geq -\frac{1}{x}\psi(a,c,x)\psi(a,c-1,x),$$
$$\psi^2(a,c,x)-\psi(a,c-1,x)\psi(a,c+1,x)\geq -\frac{1}{x}\psi^2(a,c,x)\psi(a+1,c+1,x),$$
where $a>0,$ $c<a+2,$ $x>0$ in the first inequality, and $a>1,$ $c< a+1,$ $x>0$ in the second inequality.

Next, using some ideas from \cite[Theorem 2]{baricz1} and the Tur\'an type inequalities for confluent hypergeometric functions of the second kind \cite{bais}, we deduce some new inequalities for confluent hypergeometric functions of the second kind.

\begin{theorem}\label{theorem4}
If $a>0>c$ and $x>0$, then the following inequalities are valid:
\begin{equation}\label{ineq1}
\left(\frac{\Gamma(a-c+1)}{\Gamma(-c)}\psi(a+1,c+1,x)\right)^{\frac{1}{a+1}}
<\left(\frac{\Gamma(a-c+1)}{\Gamma(1-c)}\psi(a,c,x)\right)^{\frac{1}{a}},
\end{equation}
and
\begin{equation}\label{ineq2}
2<\frac{\psi(a,c,x)}{\psi(a+1,c+1,x)}-\frac{1}{c}
\left(\frac{\Gamma(a-c+1)}{\Gamma(1-c)}\psi(a,c,x)\right)^{\frac{1}{a}}.
\end{equation}

Moreover, the following inequality
\begin{equation}\label{ineq3}
\left(\frac{\Gamma(a-c+1)}{\Gamma(1-c)}\psi(a,c,x)\right)^{\frac{c}{a(c+1)}}
<\left(\frac{\Gamma(a-c+1)}{\Gamma(-c)}\psi(a+1,c+1,x)\right)^{\frac{1}{a+1}}
\end{equation}
is valid for $a>0,c<-1$, $x>0$ and the inequality
\begin{equation}\label{ineq4}
\psi(a+1,c+1,x)<-\frac{1}{c}\psi(a,c,x)
\end{equation}
is valid for $a>0>c$ and $x>0$.
\end{theorem}

\section{\bf Proofs of main results}\label{proof}
\setcounter{equation}{0}

In this section we prove our main results and we also give alternative proofs for the Tur\'an type inequalities $(3.4)$, $(3.11)$, $(3.15)$ and left-hand side of $(3.16)$ in \cite{bais}.

\begin{proof}[\bf Proof of Theorem \ref{theorem2*}]
Taking into account the proof of \cite[Theorem 2]{bais}, we know that for $a>0$, $c<1$ and $x>0$ we have
\begin{equation}\label{int1}
\frac{_\psi\Delta_{a,c}(x)}{\psi^2(a,c,x)}=-\int_0^{\infty}\frac{t\varphi_{a,c}(t)}{(x+t)^2}dt,
\end{equation}
where
$$\varphi_{a,c}(t)=\frac{t^{-c}e^{-t}|\psi(a,c,te^{i\pi})|^{-2}}{\Gamma(a+1)\Gamma(a-c+1)},$$
and
$$_\psi\Delta_{a,c}(x)=\psi^2(a,c,x)-\psi(a-1,c-1,x)\psi(a+1,c+1,x).$$
Observe that
$$\left(\frac{x^2{}_\psi\Delta_{a,c}(x)}{\psi^2(a,c,x)}\right)'=
-\int_0^{\infty}\frac{2xt^2\varphi_{a,c}(t)}{(x+t)^3}dt<0,$$
for all $a>0$, $c<1$ and $x>0$. Thus, the function
$$x\mapsto \frac{x^2{}_\psi\Delta_{a,c}(x)}{\psi^2(a,c,x)}$$
is strictly decreasing and maps $(0,\infty)$ into $(-\infty,0).$ Consequently we have for all $a>0$, $c<1$ and $x>0$,
$$\zeta_{a,c}:=\lim_{x\rightarrow \infty}\frac{x^2{}_\psi\Delta_{a,c}(x)}{\psi^2(a,c,x)}
<\frac{x^2{}_\psi\Delta_{a,c}(x)}{\psi^2(a,c,x)}<\lim_{x\rightarrow \infty}\frac{x^2{}_\psi\Delta_{a,c}(x)}{\psi^2(a,c,x)}=:\eta_{a,c}.$$
Now to compute $\zeta_{a,c}$, we use the asymptotic expansion \cite[p. 508]{abramowitz}
$$\psi(a,c,x)\sim x^{-a}\left(1+a(c-a-1)\frac{1}{x}+\frac{1}{2}a(a+1)(a+1-c)(a+2-c)\frac{1}{x^2}+\cdots\right),$$
which is valid for large real $x$ and fixed $a$ and $c$. Thus, we have
\begin{eqnarray}
\frac{x^2{}_\psi\Delta_{a,c}(x)}{\psi^2(a,c,x)} &\sim & \frac{x^2}{\left(1+\alpha_1(a,c)\frac{1}{x}+\alpha_2(a,c)\frac{1}{x^2}+
\cdots\right)^2}\cdot\left(\left(1+\alpha_1(a,c)\frac{1}{x}+\alpha_2(a,c)
\frac{1}{x^2}+\cdots\right)^2\right. \nonumber \\ && -\left(1+\alpha_1(a-1,c-1)\frac{1}{x}+\alpha_2(a-1,c-1)\frac{1}{x^2}+
\cdots\right) \times \nonumber \\ && \hskip 0.25cm \left.\left(1+\alpha_1(a+1,c+1)\frac{1}{x}+\alpha_2(a+1,c+1)
\frac{1}{x^2}+\cdots\right) \right), \nonumber
\end{eqnarray}
where
$$\alpha_1(a,c)=a(c-a-1) ~~\mbox{and}~~ \alpha_2(a,c)=\frac{1}{2}a(a+1)(a+1-c)(a+2-c).$$
Since $2\alpha_1(a,c)=\alpha_1(a-1,c-1)+\alpha_1(a+1,c+1)$, it follows that in the above asymptotic expansion inside of square brackets
the factor ${1}/{x}$ vanishes and thus we just need to look at the coefficient of the factor ${1}/{x^2}$. But this is
$$2\alpha_2(a,c)+\alpha_1^2(a,c)-\alpha_2(a-1,c-1)
-\alpha_2(a+1,c+1)-\alpha_1(a-1,c-1)\alpha_1(a+1,c+1).$$
After some algebra we obtain that the above expression becomes $c-a-1$,
and thus we have that
\begin{equation}\label{lim1}
\zeta_{a,c}=\lim_{x\rightarrow \infty}\frac{x^2{}_\psi\Delta_{a,c}(x)}{\psi^2(a,c,x)}=c-a-1,
\end{equation}
that is, indeed the inequality \eqref{te1} is valid.

Note that in view of the proof of \cite[Theorem 2]{bais} we know that
\begin{equation}\label{lim1*}
\lim_{x\rightarrow 0}\frac{_\psi\Delta_{a,c}(x)}{\psi^2(a,c,x)}=\frac{1}{c}
\end{equation}
for $a>0>c$, and thus it follows that
$$\eta_{a,c}=\lim_{x\rightarrow 0}\frac{x^2._\psi\Delta_{a,c}(x)}{\psi^2(a,c,x)}=0.$$
In this case we obtain the Tur\'an type inequality \eqref{te8}.

Alternatively, \eqref{te1} can be proved as follows. Using \eqref{lim1} and \eqref{int1} we have
\begin{equation}\label{int2}
\int_0^{\infty}t\varphi_{a,c}(t)=1+a-c.
\end{equation}
Now from \eqref{int1} and \eqref{int2}, for all $a>0$, $c<1$ and $x>0$ we have
$$\frac{_\psi\Delta_{a,c}(x)}{\psi^2(a,c,x)}=-\int_0^{\infty}\frac{t\varphi_{a,c}(t)}{(x+t)^2}dt>
-\int_0^{\infty}\frac{t\varphi_{a,c}(t)}{x^2}=\frac{c-a-1}{x^2}.$$
This gives another proof of \eqref{te1}.

Moreover, again using \eqref{int1}, \eqref{int2} and the inequality $(x+t)^3>x^3$, we can get another proof of \eqref{te1}. Namely, we have that
$$\left(\frac{_\psi\Delta_{a,c}(x)}{\psi^2(a,c,x)}\right)'
=\int_0^{\infty}\frac{2t\varphi_{a,c}(t)}{(x+t)^3}dt
<\int_0^{\infty}\frac{2t\varphi_{a,c}(t)}{x^3}dt=\left(\frac{c-a-1}{x^2}\right)',$$
that is, the function
$$x\mapsto \frac{_\psi\Delta_{a,c}(x)}{\psi^2(a,c,x)}-\frac{c-a-1}{x^2}$$
is strictly decreasing on $(0,\infty)$. Since (see \cite{bais})
\begin{equation}\label{lim2}
\lim_{x\rightarrow \infty}\frac{_\psi\Delta_{a,c}(x)}{\psi^2(a,c,x)}=0,
\end{equation}
it results that
$$\frac{_\psi\Delta_{a,c}(x)}{\psi^2(a,c,x)}>\frac{c-a-1}{x^2},$$
which is exactly \eqref{te1}.

By using a similar trick as above we have that
\begin{equation}\label{int3}
\int_0^{\infty}\varphi_{a,c}(t)=1,
\end{equation}
since (see proof of \cite[Theorem 3]{bais})
$$\lim_{x\rightarrow \infty}\frac{_\psi\Delta_a(x)}{\psi^2(a,c,x)}=0,$$
and
$$\frac{(1+a-c)_\psi\Delta_a(x)}{\psi^2(a,c,x)}=1-\int_0^{\infty}\frac{x^2\varphi_{a,c}(t)}{(x+t)^2}dt,$$
where $_\psi\Delta_a(x)=\psi^2(a,c,x)-\psi(a-1,c,x)\psi(a+1,c,x)$, $a>0$, $c<1$ and $x>0$.

Now in view of \eqref{int3} and the inequality $2tx<(t+x)^2$ we obtain
$$\frac{_\psi\Delta_{a,c}(x)}{\psi^2(a,c,x)}=-\int_0^{\infty}\frac{t\varphi_{a,c}(t)}{(x+t)^2}dt > -\frac{1}{2x}\int_0^{\infty}\varphi_{a,c}(t)dt=-\frac{1}{2x};$$
or equivalently
$$\psi^2(a,c,x)-\psi(a-1,c-1,x)\psi(a+1,c+1,x)>-\frac{1}{2x}\psi^2(a,c,x),$$
where $a>0$, $c<1$ and $x>0$. This proves the Tur\'an type inequality \eqref{te2}.

Taking into account of the differentiation formula \cite[p. 507]{abramowitz}
\begin{equation}\label{derivative}
\psi'_{a,c,x}(x)=-a\psi(a+1,c+1,x),
\end{equation}
we obtain
\begin{align*}
&\left(\frac{_\psi\Delta_{a,c}(x)}{\psi^2(a,c,x)}\right)'\\
&\ \ \ \ =\frac{(a-1)\psi^2(a,c,x)\psi(a+1,c+1,x)
+(a+1)\psi(a-1,c-1,x)\psi(a+2,c+2,x)\psi(a,c,x)}{\psi^3(a,c,x)}\\
&\ \ \ \ -\frac{2a\psi(a-1,c-1,x)\psi^2(a+1,c+1,x)}{\psi^3(a,c,x)}
\end{align*}
and in view of the asymptotic expansion \cite[p. 508]{abramowitz}
\begin{equation}\label{asymp}
\psi(a,c,x)\sim \frac{\Gamma(1-c)}{\Gamma(a-c+1)},
\end{equation}
where $a>0$ and $c<1$ are fixed and $x\rightarrow 0$, we obtain that
$$\left(\frac{_\psi\Delta_{a,c}(x)}{\psi^2(a,c,x)}\right)'\sim \frac{(a-1)\Gamma^2(1-c)\Gamma(-c)+(a+1)\Gamma(1-c)\Gamma(2-c)\Gamma(-1-c)-2a\Gamma(2-c)\Gamma^2(-c)}
{\Gamma^3(1-c)},$$
where $a>1$ and $c<-1$ are fixed and $x\rightarrow 0$. After some algebra (using $\Gamma(\alpha+1)=\alpha\Gamma(\alpha)$) we get
$$\left[\frac{_\psi\Delta_{a,c}(x)}{\psi^2(a,c,x)}\right]'\sim \frac{2(c-a)}{c^2(c+1)},$$
where $a>1$ and $c<-1$. On the other hand, using \eqref{int1} we obtain
$$\lim_{x\rightarrow 0}\left(\frac{_\psi\Delta_{a,c}(x)}{\psi^2(a,c,x)}\right)'
=2\int_0^{\infty}\frac{\varphi_{a,c}(t)}{t^2}dt.$$
Thus we have
\begin{equation}\label{int4}
\int_0^{\infty}\frac{\varphi_{a,c}(t)}{t^2}dt=\frac{c-a}{c^2(c+1)},
\end{equation}
where $a>1$ and $c<-1$.
Now in view of \eqref{int4} and the inequality $(x+t)^3>t^3$ we have
$$\left(\frac{_\psi\Delta_{a,c}(x)}{\psi^2(a,c,x)}\right)'=
\int_0^{\infty}\frac{2t\varphi_{a,c}(t)}{(x+t)^3}dt<
2\int_0^{\infty}\frac{\varphi_{a,c}(t)}{t^2}dt=
\left(\frac{2(c-a)}{c^2(c+1)}x\right)',$$
and thus the function
$$x\mapsto \frac{_\psi\Delta_{a,c}(x)}{\psi^2(a,c,x)}-2x\left(\frac{c-a}{c^2(c+1)}\right)$$
is strictly decreasing on $(0,\infty)$ for $a>1$ and $c<-1$. Consequently in view of \eqref{lim1*} the Tur\'an type inequality \eqref{te1r} follows.

Finally, observe that in view of \eqref{lim1*} the inequality \eqref{te1r} is sharp as $x\rightarrow 0,$ while \eqref{lim2} gives the sharpness of \eqref{te1} and \eqref{te2} as $x\rightarrow \infty$.
\end{proof}

\begin{proof}[\bf Proof of Theorem \ref{theorem3*}]

Let us consider the Tur\'anian
$$_\psi\Delta_a(x)=\psi^2(a,c,x)-\psi(a-1,c,x)\psi(a+1,c,x).$$
According to the proof of \cite[Theorem 3]{bais} we know that
\begin{equation}\label{int5}
\frac{(1+a-c)_\psi\Delta_a(x)}{\psi^2(a,c,x)}=
1-\int_0^{\infty}\frac{x^2\varphi_{a,c}(t)}{(x+t)^2}dt,
\end{equation}
where $a>0$, $c<1$ and $x>0$.

Now consider the function $x\mapsto x^2\left(\frac{(1+a-c)_\psi\Delta_a(x)}{\psi^2(a,c,x)}-1\right)$. Since for all $x>0$ and $t>0$ we have
$$\left(\frac{x^4}{(x+t)^2}\right)'=\frac{2x^3(x+2t)}{(x+t)^3}>0,$$
and thus the above function is strictly increasing on $(0,\infty)$.
Since (see the proof of \cite[Theorem 3]{bais})
$$\lim_{x\rightarrow 0}\frac{(1+a-c){}_\psi\Delta_a(x)}{\psi^2(a,c,x)}=1,$$
for $a>1>c$ and $x>0$, it follows that
$$\lim_{x\rightarrow 0}x^2\left(\frac{(1+a-c){}_\psi\Delta_a(x)}{\psi^2(a,c,x)}-1\right)=0,$$
and thus we have the inequality \eqref{te10}, which is exactly the left-hand side of (3.15) in \cite{bais}.

By using \eqref{int2}, \eqref{int5} and the inequality $\frac{2x}{(x+t)^3}<\frac{2}{x^2}$ we have
$$\left(\frac{_\psi\Delta_a(x)}{\psi^2(a,c,x)}\right)'=
-\frac{1}{1+a-c}\int_0^{\infty}\frac{2xt\varphi_{a,c}(t)}{(x+t)^3}dt>
-\frac{1}{1+a-c}\int_0^{\infty}\frac{2t\varphi_{a,c}(t)}{x^2}dt=\left(\frac{2}{x}\right)',$$
that is, the function
$$x\mapsto \frac{\psi\Delta_a(x)}{\psi^2(a,c,x)}-\frac{2}{x}$$
is strictly increasing on $(0,\infty)$. Since (see the proof of \cite[Theorem 3]{bais})
$$\lim_{x\rightarrow \infty}\frac{\psi\Delta_a(x)}{\psi^2(a,c,x)}=0,$$
we deduce that
$$\psi^2(a,c,x)-\psi(a-1,c,x)\psi(a+1,c,x)<\frac{2}{x}\psi^2(a,c,x),$$
where $a>0$, $c<1$ and $x>0$. This proves the right-hand side of \eqref{te3}.

We also note that \eqref{int3} can be used to get another proof of \eqref{te9}, which is the right-hand side of (3.15) in \cite{bais}. Namely, for $a>0$, $c<1$ and $x>0$ we have
$$\frac{(1+a-c)_\psi\Delta_a(x)}{\psi^2(a,c,x)}=
1-\int_0^{\infty}\frac{x^2\varphi_{a,c}(t)}{(x+t)^2}dt>1-\int_0^{\infty}\varphi_{a,c}(t)=0,$$
which gives right-hand side of $(3.15)$ in \cite{bais}.

Now, observe that in view of \eqref{int1} and \eqref{lim1*} we have
\begin{equation}\label{int6}
\lim_{x\rightarrow 0}\frac{_\psi\Delta_{a,c}(x)}{\psi^2(a,c,x)}=\frac{1}{c}=-\int_0^{\infty}\frac{\varphi_{a,c}(t)}{t}dt,
\end{equation}
where $a>0>c$. From \eqref{int6} and the inequality $(x+t)^2>2xt$
we have for all $a>0>c$ and $x>0$,
$$\frac{(1+a-c)_\psi\Delta_a(x)}{\psi^2(a,c,x)}=
1-\int_0^{\infty}\frac{x^2\varphi_{a,c}(t)}{(x+t)^2}dt>
1-\frac{x}{2}\int_0^{\infty}\frac{\varphi_{a,c}(t)}{t}dt=1+\frac{x}{2c},$$
which proves the left-hand side of \eqref{te3}.

Now using \eqref{int5}, inequality $(x+t)^3>t^3$ and \eqref{int4}, we have for all $a>1$ and $c<-1$
$$\left(\frac{(1+a-c)_\psi\Delta_a(x)}{\psi^2(a,c,x)}\right)'=
-\int_0^{\infty}\frac{2xt\varphi_{a,c}(t)}{(x+t)^3}dt>
-2x\int_0^{\infty}\frac{\varphi_{a,c}(t)}{t^2}dt=
\left(-\left(\frac{c-a}{c^2(c+1)}\right)x^2\right)'.$$
Thus, the function
$$x\mapsto \frac{(1+a-c)_\psi\Delta_a(x)}{\psi^2(a,c,x)}+\left(\frac{c-a}{c^2(c+1)}\right)x^2$$
is strictly increasing on $(0,\infty)$ and hence using the limit
$$\lim_{x\rightarrow 0} \frac{(1+a-c)_\psi\Delta_a(x)}{\psi^2(a,c,x)}=1,$$
the Tur\'an type inequality \eqref{te5} follows.

Since (see proof of \cite[Theorem 3]{bais})
$$\lim_{x\rightarrow 0}\frac{_\psi\Delta_a(x)}{\psi^2(a,c,x)}=\frac{1}{1+a-c}$$
and
$$\lim_{x\rightarrow 0}\frac{_\psi\Delta_a(x)}{\psi^2(a,c,x)}=0,$$
the left-hand side of \eqref{te3}, and the inequality \eqref{te5} are clearly sharp as $x\rightarrow 0,$ while the right-hand side of \eqref{te3} is sharp as $x\rightarrow \infty$.
\end{proof}

\begin{proof}[\bf Proof of Theorem \ref{theorem4*}]
By using \eqref{int2}, the inequality $(x+t)^3>x^3$ and the next relation (see the proof of \cite[Theorem 4]{bais})
$$\left(\frac{_\psi\Delta_c(x)}{\psi^2(a,c,x)}\right)'=
\frac{a}{1+a-c}\int_0^{\infty}\frac{2t\varphi_{a,c}(t)}{(x+t)^3}dt,$$
where
$${}_\psi\Delta_c(x)=\psi^2(a,c,x)-\psi(a,c-1,x)\psi(a,c+1,x),$$
we get for all $a>0$, $c<1$ and $x>0$
$$\left(\frac{{}_\psi\Delta_c(x)}{\psi^2(a,c,x)}\right)'<\left(-\frac{a}{x^2}\right)',$$
that is, the function
$$x\mapsto \frac{{}_\psi\Delta_c(x)}{\psi^2(a,c,x)}+\frac{a}{x^2}$$
is strictly decreasing on $(0,\infty)$. Since (see the proof of \cite[Theorem 4]{bais})
\begin{equation}\label{lim3}
\lim_{x\rightarrow \infty}\frac{{}_\psi\Delta_c(x)}{\psi^2(a,c,x)}=0
\end{equation}
and hence we get the inequality \eqref{te6}.

By using \eqref{int6} and the inequality $(x+t)^2>t^2$ we obtain
$$\frac{{}_\psi\Delta_c(x)}{\psi^2(a,c,x)}=
-\frac{a}{1+a-c}\int_0^{\infty}\frac{t\varphi_{a,c}(t)}{(x+t)^2}dt>
-\frac{a}{1+a-c}\int_0^{\infty}\frac{\varphi_{a,c}(t)}{t}dt=\frac{a}{c(1+a-c)},$$
where $a>0>c$ and $x>0$. This is the left-hand side of (3.16) in \cite{bais}, that is, the inequality \eqref{te11}.

Finally, using \eqref{int4} for $a>1$, $c<-1$ and $x>0$ we have
$$\left(\frac{_\psi\Delta_c(x)}{\psi^2(a,c,x)}\right)'=
\frac{a}{1+a-c}\int_0^{\infty}\frac{2t\varphi_{a,c}(t)}{(x+t)^3}dt<
\frac{2a}{1+a-c}\int_0^{\infty}\frac{\varphi_{a,c}(t)}{t^2}dt=
\left(\frac{2a}{1+a-c}\left(\frac{c-a}{c^2(c+1)}\right)x\right)',$$
that is, the function
$$x\mapsto \frac{_\psi\Delta_c(x)}{\psi^2(a,c,x)}-\frac{2a}{1+a-c}\left(\frac{c-a}{c^2(c+1)}\right)x$$
is decreasing on $(0,\infty)$. Since (see the proof of \cite[Theorem 4]{bais})
\begin{equation}\label{lim4}
\lim_{x\rightarrow 0}\frac{_\psi\Delta_c(x)}{\psi^2(a,c,x)}=\frac{a}{c(1+a-c)},
\end{equation}
we get the inequality \eqref{te6r}.

Now in view of \eqref{lim3}, the inequality \eqref{te6} is sharp as $x\rightarrow \infty$ while \eqref{lim4} gives the sharpness of \eqref{te6r} as $x\rightarrow 0$.
\end{proof}

\begin{proof}[\bf Proof of Theorem \ref{theorem4}]
To prove \eqref{ineq1}, consider the function $f_{a,c}:(0,\infty)\rightarrow \mathbb{R}$, defined by
$$f_{a,c}(x)=\frac{1}{a}\log \psi(a,c,x)-\frac{1}{a+1}\log \psi(a+1,c+1,x),$$
which in view of the Tur\'an type inequality \eqref{te8} and the differentiation formula \eqref{derivative} gives
$$f_{a,c}'(x)=\frac{\psi(a+2,c+2,x)}{\psi(a+1,c+1,x)}-\frac{\psi(a+1,c+1,x)}{\psi(a,c,x)}>0.$$
Therefore $f_{a,c}$ is increasing on $(0,\infty)$ and hence $f_{a,c}(x)>f_{a,c}(0)$ for $x>0$ and consequently by using the asymptotic expansion \eqref{asymp}, the inequality \eqref{ineq1} follows.

The inequality \eqref{ineq2} follows from \eqref{ineq1} and arithmetic-geometric mean inequality. Namely, by using inequality \eqref{ineq1} twice, we have
\begin{align*}
\frac{\psi(a+1,c+1,x)}{\psi(a,c,x)}&<
\left(\frac{\Gamma(-c)}{\Gamma(a-c+1)}\right)^{\frac{1}{a+1}}
\left(\frac{\Gamma(a-c+1)}{\Gamma(1-c)}\right)^{\frac{1}{a}}
\left(\psi(a,c,x)\right)^{\frac{1-a}{a}}\left(\psi(a+1,c+1,x)\right)^{\frac{a}{a+a}}\\
&<-\frac{1}{c}\left(\frac{\Gamma(a-c+1)}{\Gamma(1-c)}\psi(a,c,x)\right)^{\frac{1}{a}},
\end{align*}
which by using the following arithmetic-geometric mean inequality
$$2\leq \frac{\psi(a,c,x)}{\psi(a+1,c+1,x)}+\frac{\psi(a+1,c+1,x)}{\psi(a,c,x)}$$
gives the inequality \eqref{ineq2}.

Now, to prove the inequality \eqref{ineq3}, consider the function $g_{a,c}:(0,\infty)\rightarrow \mathbb{R}$, defined by
$$g_{a,c}(x)=\frac{c}{a(c+1)}\log \psi(a,c,x)-\frac{1}{(a+1)}\log \psi(a+1,c+1,x),$$
which in view of the Tur\'an type inequality \eqref{te7} and the differentiation formula \eqref{derivative} gives
$$g_{a,c}'(x)=\frac{\psi(a+2,c+2,x)}{\psi(a+1,c+1,x)}-\frac{c}{c+1}\frac{\psi(a+1,c+1,x)}{\psi(a,c,x)}<0.$$
Therefore $g_{a,c}$ is decreasing on $(0,\infty)$ and $g_{a,c}(x)<g_{a,c}(0)$ for $x>0$ and hence by \eqref{asymp} the inequality \eqref{ineq3} follows.

To prove the inequality \eqref{ineq4}, consider the function $h_{a,c}:(0,\infty)\rightarrow \mathbb{R}$, defined by
$$h_{a,c}(x)=\log \psi(a,c,x)-\log \psi(a+1,c+1,x),$$
which according to the following Tur\'an type inequality for confluent hypergeometric functions of the second kind \cite{bais}
\begin{equation}\label{te111}\psi^2(a,c,x)-\psi(a-1,c-1,x)\psi(a+1,c+1,x)<\frac{1}{a}\psi^2(a,c,x),\end{equation}
valid for $a>0$, $c\in \mathbb{R}$ and $x>0$, and the differentiation formula \eqref{derivative} gives
$$h_{a,c}'(x)=\frac{\psi(a+2,c+2,x)}{\psi(a+1,c+1,x)}-
\left(\frac{a}{a+1}\right)\frac{\psi(a+1,c+1,x)}{\psi(a,c,x)}>0.$$
Hence $h_{a,c}$ is increasing on $(0,\infty)$ and consequently in view of \eqref{asymp} inequality \eqref{ineq4} follows.

We note that the Tur\'an type inequality \eqref{te111} has been proved in \cite{bais} using the classical H\"{o}lder-Rogers inequality for integrals. Alternatively, this inequality can be proved using convolution theorem for Laplace transforms. Namely, by taking into account of the integral representation \eqref{integr} and the convolution theorem for Laplace transforms, it follows that
\begin{eqnarray*}
_\psi\Delta_{a,c}(x)&=&\psi^2(a,c,x)-\psi(a-1,c-1,x)\psi(a+1,c+1,x)\\
&=&\frac{1}{\Gamma^2(a)}\int_0^{\infty}e^{-xt}
\left(\int_0^t(t-u)^{a-1}(1+t-u)^{c-a-1}u^{a-1}(1+u)^{c-a-1}du\right)dt\\
&&-\frac{1}{\Gamma(a-1)\Gamma(a+1)}\int_0^{\infty}e^{-xt}
\left(\int_0^t(t-u)^{a-2}(1+t-u)^{c-a-1}u^a(1+u)^{c-a-1}du\right)dt\\
&=&\frac{1}{\Gamma^2(a)}\int_0^{\infty}e^{-xt}
\left(\int_0^t(t-u)^{a-2}(1+t-u)^{c-a-1}u^{a-1}(1+u)^{c-a-1}\left(t-u\frac{2a-1}{a}\right)du\right)dt.
\end{eqnarray*}
Now by using the change of variable $u=\frac{t(1+y)}{2}$, it follows that,
\begin{eqnarray*}
&&\int_0^t(t-u)^{a-2}(1+t-u)^{c-a-1}u^{a-1}(1+u)^{c-a-1}\left(t-u\frac{2a-1}{a}\right)du
\\&&=\int_{-1}^1\left(\frac{t(1-y)}{2}\right)^{a-2}
\left(1+\frac{t(1-y)}{2}\right)^{c-a-1}
\left(\frac{t(1+y)}{2}\right)^{a-1}
\left(1+\frac{t(1+y)}{2}\right)^{c-a-1}\frac{y(1-2a)+1}{2a}\frac{t^2}{2}dy\\
&&=2\int_0^1\frac{t^{2a-1}}{2^{2a-1}}
\left(1-y^2\right)^{a-2}\left(\left(1+\frac{t}{2}\right)^2-
\frac{t^2y^2)}{4}\right)^{c-a-1}\frac{1}{a}\left(1+(1-2a)y^2\right)dy.
\end{eqnarray*}
On the other hand, by using similar steps, $\psi^2(a,c,x)/a$ can be rewritten as
\begin{align*}
&\frac{1}{a\Gamma^2(a)}\int_0^{\infty}e^{-xt}
\left(\int_0^t(t-u)^{a-1}(1+t-u)^{c-a-1}u^{a-1}(1+u)^{c-a-1}du\right)dt\\
&=\frac{1}{a\Gamma^2(a)}\int_0^{\infty}e^{-xt}\int_{-1}^1
\left(\frac{t^2(1-y)(1+y)}{4}\right)^{a-1}\left(1+\frac{t(1-y)}{2}\right)^{c-a-1}\left(1+\frac{t(1+y)}{2}\right)^{c-a-1}\frac{t}{2}dydt\\
&=\frac{1}{a\Gamma^2(a)}\int_0^{\infty}e^{-xt}
\left(2\int_0^1\frac{t^{2a-1}}{2^{2a-1}}\left(1-y^2\right)^{a-2}
\left(\left(1+\frac{t}{2}\right)^2-\frac{t^2y^2)}{4}\right)^{c-a-1}\frac{1}{a}
\left(1-y^2\right)dy\right)dt.
\end{align*}
Since for all $a>1$ and $y\in(0,1)$, we have
$$\frac{1}{a}\left(1+(1-2a)y^2\right)<\frac{1}{a}(1-y^2),$$
it follows that $_\psi\Delta_{a,c}(x)<\frac{1}{a}\psi^2(a,c,x)$, that is, we have the
Tur\'an type inequality \eqref{te111} for $a>1$, $c\in \mathbb{R}$ and $x>0,$ which is exactly the inequality (3.11) in \cite{bais}.
\end{proof}

\end{document}